\documentstyle[12pt]{article}
\title{On the cohomology of $s\ell(m+1,\rit)$ acting on differential operators and
$s\ell(m+1,\rit)$-equivariant symbol}
\author{P.B.A. Lecomte }
\date{}
\newtheorem{theo}{\indent Theorem}[section]
\newtheorem{prop}[theo]{\indent Proposition}
\newtheorem{lema}[theo]{\indent Lemma}
\newtheorem{coro}[theo]{\indent Corollary}
\newtheorem{rema}[theo]{\indent Remark}

\def\mod#1{\left|{#1}\right|}
\newcommand{\rit}{{\rm I\!R}}
\newcommand{\nit}{{\rm I\!N}}
\def\A{{\cal A}}
\def\B{{\cal B}}
\def\C{{\cal C}}
\def\D{{\cal D}}
\def\S{{\cal S}}
\def\L{{\cal L}}
\def\T{{\cal T}}
\newcommand{\cqfd}{\hspace*{\fill}\rule{2mm}{2mm}} 
\def\im{\mathop{\rm im}}

\def\scal#1#2{\left\langle{#1},{#2}\right\rangle}

\def\resume{\if@twocolumn
\section*{RŽsumŽ}
\else \small 
\begin{center}
{\bf 
RŽsumŽ
\vspace{-.5em}\vspace{0pt}} 
\end{center}
\quotation 
\fi}
\def\endresume{\if@twocolumn\else\endquotation\fi}

\begin{document}
\maketitle
\begin{abstract}
Let $\D_{\lambda\mu}^k(\rit^m)$ denote the space of differential operators of order $\leq
k$ from the space of $\lambda$-densities of $\rit^m$ into that of $\mu$-densities and let
$\S_\delta^k(\rit^m)$ be the space of $k$-contravariant, symmetric tensor fields on
$\rit^m$ valued in the $\delta$-densities, $\delta = \mu-\lambda$.

Denote by $s\ell_{m+1}$ the projective embedding of $s\ell(m+1,\rit)$ as a subalgebra of
the Lie algebra of vector fields of $\rit^m$.

One computes the cohomology of $s\ell_{m+1}$ with coefficients in the space of
differential operators from $\S_\delta^p(\rit^m)$ into $\S_\delta^q(\rit^m)$. It is non
vanishing only for some {\it critical} values of $\delta$. For $m=1$, these are the values
pointed out by H. Gargoubi in a completely different context \cite{3}.

This allows to determine the condition under which the short exact sequence of
$s\ell_{m+1}$-modules
$$0 \rightarrow \D_{\lambda\mu}^{k-1}(\rit^m) \rightarrow \D_{\lambda\mu}^k(\rit^m)
\stackrel{\sigma}{\rightarrow} \S_\delta^k(\rit^m) \rightarrow 0$$
is split ($\sigma$ is the symbol map).

From this one recovers and generalizes useful results about the structure of the
$s\ell_{m+1}$-module $\D_{\lambda\mu} (\rit^m) = \oplus_k \D_{\lambda\mu}^k(\rit^m)$
\cite{1,2,3,4}.

The cohomology of the latter is also computed.
\end{abstract}

\noindent {\it Key words} : Lie algebra, differential operators, cohomology, infinitesimal
homographies.

\medskip\noindent
Classification numbers : 58B99, 17B56, 17B65, 47E05.

\medskip\noindent
Titre courant : cohomology of $s\ell(m+1,\rit)$ on differential operators.

\begin{resume}
Soient $\D_{\lambda\mu}^k(\rit^m)$ l'espace des opŽrateurs diffŽrentiels d'ordre $\leq k$
sur les $\lambda$-densitŽs de $\rit^m$ ˆ valeurs dans les $\mu$-densitŽs de $\rit^m$ et
$\S_\delta^k(\rit^m)$ l'espace des champs de tenseurs $k$-contravariants symŽtriques sur
$\rit^m$ ˆ valeurs dans les $\delta$-densitŽs, o $\delta = \mu-\lambda$.

Soit Žgalement $s\ell_{m+1}$ le plongement projectif de $s\ell(m+1,\rit)$ dans l'algbre
de Lie des champs de vecteurs de $\rit^m$.

On calcule la cohomologie de $s\ell_{m+1}$ ˆ coefficients dans les opŽrateurs
diffŽrentiels de $\S_\delta^p(\rit^m)$ dans $\S_\delta^q(\rit^m)$. Elle n'est non 
nulle que pour certaines valeurs, dites {\it critiques}, de $\delta$. Pour $m=1$, ce sont
les valeurs mises en Žvidence par H. Gargoubi dans un tout autre contexte \cite{3}.

A l'aide des rŽsultats obtenus, on dŽtermine ˆ quelle condition la courte suite exacte de
$s\ell_{m+1}$-modules
$$0 \rightarrow \D_{\lambda\mu}^{k-1}(\rit^m) \rightarrow \D_{\lambda\mu}^k(\rit^m)
\stackrel{\sigma}{\rightarrow} \S_\delta^k(\rit^m) \rightarrow 0$$
est scindŽe ($\sigma$ est le passage au symbole).

Cela permet de retrouver et de gŽnŽraliser diffŽrents rŽsultats utiles concernant la
structure de $s\ell_{m+1}$-module de $\D_{\lambda\mu}(\rit^m) = \oplus_k
\D_{\lambda\mu}^k(\rit^m)$ \cite{1,2,3,4}.

La cohomologie de ce dernier est Žgalement calculŽe.

\end{resume}

\section{Introduction}
\indent\indent Let $\D_{\lambda \mu}(M)$ denote the space of differential operators from
the space of scalar $\lambda$-densities over a smooth manifold $M$ into the space of
scalar $\mu$-densities.

The Lie derivation with respect to vector fields equips $\D_{\lambda\mu}(M)$ with a
natural structure of module over the Lie algebra of vector fields $Vect(M)$ of $M$. The
study of that module has been started in \cite{1,2,3,4,5}. See also \cite{6} for
pseudodifferential operators on the complex line. As a vector space, $\D_{\lambda\mu}(M)$
is isomorphic to the graded space
$$\S_\delta(M) = \bigoplus_{k \geq 1}
\frac{\D_{\lambda\mu}^k(M)}{\D_{\lambda\mu}^{k-1}(M)}$$
associated to the filtration of $\D_{\lambda\mu}(M)$ by the order $k$ of differentiation.
The $k$-th term of this sum is nothing else but the space $\S_\delta^k(M)=\Gamma(\vee^kTM
\otimes \Delta^\delta M)$ of $k$-symmetric contravariant tensor fields valued in the scalar
densities of weight $\delta=\mu-\lambda$ ($\Delta^\delta M$ denotes the vector bundle of
$\delta$-densities of $M$). The isomorphism between
$\D_{\lambda\mu}^k(M)/\D_{\lambda\mu}^{k-1}(M)$ and $\S_\delta^k$ is of course induced
in a natural way by the symbol map
$$\sigma : \D_{\lambda\mu}^k(M) \rightarrow \S_\delta^k(M).$$

This map commutes with the Lie derivation of operators on the left and of tensor fields on
the right. But, however, $\D_{\lambda\mu}(M)$ and $\S_\delta(M)$ are not isomorphic as
$Vect(M)$-modules~: $\D_{\lambda\mu}(M)$ is better viewed as a non trivial deformation of
$\S_\delta(M)$.

For $\lambda=\mu$ and $m>1$, it has been shown \cite{4}, nevertheless, that
$\D_{\lambda\mu}(M)(\rit^m)$ and $\S_0(\rit^m)$ are isomorphic {\it as
$s\ell_{m+1}$-modules}, where $s\ell_{m+1} \subset Vect(\rit^m)$ is the canonical embedding
of $s\ell(m+1,\rit)$ as the algebra of infinitesimal linear fractional local
transformations of $\rit^m$ (note that it is a maximal subalgebra of the algebra of
polynomial vector fields on $\rit^m$).

A similar result has been obtained in the one dimensional case (cf. \cite{6} for
pseudodifferential operators, \cite{3} for differential operators)~: except for special
critical values of $\lambda$ and $\mu$, $\D_{\lambda\mu}(M)$ is again
$s\ell_2$-isomorphic to $\S_\delta(M)$. For the critical values, $\D_{\lambda\mu}(M)$ is
a non trivial deformation of $\S_\delta(M)$ and it is described in \cite{3} in terms of a
family of cocycles, taylored on the occasion, of $s\ell_2$ acting on the space of
differential operators from $\Gamma(\Delta^{-n/2}M)$ into $\Gamma(\Delta^{(n+2)/2}M)$.

In both cases, $m=1$ and $m>1$, the $s\ell_{m+1}$-module structure of
$\D_{\lambda\mu}(\rit^m)$ proved to be a powerful tool in studying its
$Vect(\rit^m)$-structure.

In some sense, the possible difference between the $L$-modules $\D_{\lambda\mu}(M)$ and
$\S_\delta(M)$, $L=Vect(M)$ or $s\ell_{m+1}$, follows from the fact that the short exact
sequence of modules
\begin{equation} \label{E0}
0 \rightarrow \D_{\lambda\mu}^{k-1}(M) \stackrel{i}{\rightarrow}
\D_{\lambda\mu}^{k}(M) \stackrel{\sigma}{\rightarrow} \S_\delta^k(M) \rightarrow 0
\end{equation}
($i$ : the inclusion) maybe is not split. The obstruction against splitting is a member of
$$H^1 (L, Hom(\S_\delta^k(M), \D_{\lambda\mu}^{k-1}(M))).$$
Composing with the symbol map $\sigma$ : $\D_{\lambda\mu}^{k-1}(M) \rightarrow
\S_{\delta}^{k-1}(M)$, one gets a sort of first order approximation of this
obstruction, namely an element of
$$H^1(L,Hom(\S_\delta^k(M),\S_\delta^{k-1}(M))).$$

In this paper, we compute the cohomology space
$$H(s\ell_{m+1}, \D(\S_\delta^p(\rit^m),\S_\delta^q(\rit^m))).$$

This allows us not only to recover the critical values of \cite{3} but also to find
critical values for $m>1$. 
To do that, we compute the spaces
$$H(s\ell_{m+1},\D^k(\S_\delta^p(\rit^m),\S_\delta^q(\rit^m)))$$
$(\D^k(A,B)$ denotes the space of differential operators of order $\leq k$ from $A$ into
$B$).

They are computed by induction on $k$,
using the short exact sequence associated to the inclusion of the operators of order $k-1$
into these of order $k$. The critical values occur when $p-q \in \nit_0$. They are then
these for which the connecting homomorphism of the corresponding exact sequence in
cohomology vanishes. Equivalently, they are these for which
$$H(s\ell_{m+1},\D(\S_\delta^p(\rit^m),\S_\delta^q(\rit^m))) \neq 0.$$

We also compute the space $Hom_{s\ell_{m+1}} (\S_\delta^p(\rit^m),\S_\delta^q(\rit^m))$.
This is quite hard because one is left to show that a mapping from $\S_\delta^p(\rit^m)$
into $\S_\delta^q(\rit^m)$ is local knowing only that it is commuting with the $L_X$, $X
\in s\ell_{m+1}$. As a corollary, we get  generalization of a result of \cite{4}~:
if $(m+1)\delta-m \not\in \nit_0$, then a linear map from $\S_\delta(\rit^m)$ is
$s\ell_{m+1}$-equivariant if and only if its restriction to each $\S_\delta^k(\rit^m)$, $k
\in \nit$, is a constant multiple of the identity.

As a consequence of our computations, we find the necessary and sufficient condition for
(\ref{E0}) to be split as a sequence of $s\ell_{m+1}$-modules (with $M=\rit^m$). In
particular, we show that if $(m+1)\delta -m\not\in \nit_0$, then the $s\ell_{m+1}$-modules
$\D_{\lambda\mu}(\rit^m)$ and $\S_\delta(\rit^m)$ are canonically isomorphic. This also
generalizes a useful result of \cite{4}.

In fact, the above mentioned cohomology classes constructed in \cite{3} are valued in
$\D_{- \frac{n}{2}, 1+ \frac{n}{2}}(\rit)$. We recover them as a consequence of our
computations due to the fact that $\D(\S_\delta^p(\rit),\S_\delta^q(\rit))$ is isomorphic
to $\D_{\delta-p,\delta-q}(\rit^m)$.

Besides, very little more is needed then to compute the cohomology of $s\ell_{m+1}$ with
coefficients in $\D_{\lambda\mu}(\rit^m)$, $m \geq 1$. This is done in the last section.
The difference between the cases $m>1$ and $m=1$ is then really significant.

\section{The projective embedding of $s\ell(m+1,\rit)$}

\indent\indent For each $T \in s\ell (m+1,\rit)$ we denote by $T^\ast$ the fundamental
vector fields on $P_m\rit$ associated to $T$ and to the canonical action of $SL(m+1,\rit)$
by linear projective transformations. We also denote by $T^\ast$ its restriction to
$\rit^m$, viewed as the open set of points of $P_m \rit$ having a non vanishing $(m+1)$-th
projective coordinate.

The mapping $T \rightarrow T^\ast$ is an injective homomorphism of Lie algebras from
\linebreak
$s\ell(m+1,\rit)$ into the Lie algebra $Vect(\rit^m)$ of vector-fields on $\rit^m$.

Following \cite{4}, we call it the {\it projective embedding of} $s\ell(m+1,\rit)$. We now
describe it in a way well suited to our purpose.

We realize $s\ell(m+1,\rit)$ as $\rit^m \oplus g\ell(m,\rit) \oplus \rit^{m\ast}$, where
$\rit^m$ and $\rit^{m\ast}$ are abelian subalgebras, where the adjoint action of $A \in
g\ell(m,\rit)$ is the natural representation of $g\ell(m,\rit)$ and where the bracket of
$h \in \rit^m$ and $\alpha \in \rit^{m\ast}$ is given by
$$[h,\alpha] = \alpha(h){\bf 1}+h \otimes \alpha$$
($\bf 1$ is the unit matrix). It is easily seen that the fundamental vector fields
associated to $h=(h^i) \in \rit^m$, $A \in (A_j^i) \in g\ell(m,\rit)$ and $\alpha =
(\alpha_i) \in \rit^{m\ast}$ are
$$h^\ast =-h^i \partial_i, \quad A^\ast = -A_j^ix^j\partial_i\quad \mbox{and} \quad
\alpha^\ast = \alpha(x)x^i \partial_i$$
(sums over repeated indices are understood and $\partial_i$ denotes the partial derivative
with respect to $x^i$).

Note that $s\ell(m+1,\rit)$ is a graded Lie algebra, $\rit^m$, $g\ell(m,\rit)$ and
$\rit^{m\ast}$ being the homogeneous components of degree $-1,0$ and 1 respectively.

For the sake of brevity, we denote $s\ell_{m+1}$ the above realization of
$s\ell(m+1,\rit)$.

\section{Cohomology of $s\ell_{m+1}$ associated to a representation of $g\ell(m,\rit)$}

\indent\indent Let $(V,\rho)$ be a finite dimensional representation of $g\ell(m,\rit)$.

The space $C_\infty(\rit^m,V)$ of smooth $V$-valued functions on $\rit^m$ is a
representation of $Vect(\rit^m)$. The action $L_X^\rho$ of a vector field $X$ on a function
$f$ is given by
$$L_X^\rho f = X.f-\rho(DX)f$$
where $X.f$ is the usual derivative $X^i\partial_if$ and where $DX=(\partial_jX^i)$ is the
differential of $X$.

In fact, $C_\infty(\rit^m,V)$ is the space of sections of the (trivial) bundle associated
to the linear frame bundle of $\rit^m$ and to $V$ and $L^\rho$ is the corresponding natural
Lie derivation.

Our purpose is to compute the Chevalley-Eilenberg cohomology of the restriction of $L^\rho$
to $s\ell_{m+1}$.

Recall that $L^\rho$ induces a representation of $\L^\rho$ on the space of cochains. It is
defined by
$$\L_X^\rho = i_X \circ \partial^\rho + \partial^\rho \circ i_X, \quad \forall X \in
s\ell_{m+1},$$
where $\partial^\rho$ is the coboundary operator. One has also
$$(\L_X^\rho c)(X_0, \ldots, X_{s-1}) = L_X^\rho(c(X_0, \ldots, X_{s-1})) - \sum_{0 \leq i
< s} c(X_0, \ldots, [X,X_i], \ldots , X_{s-1})$$
for each $s$-cochain $c$.

\begin{prop} \label{P2.1}
For each $u \in \{0, \ldots,m\}$, there exists a unique linear 
\linebreak 
$\chi^u :\Lambda(g\ell(m,\rit), \Lambda^u(\rit^{m,\ast},V)) \rightarrow
\Lambda(s\ell_{m+1}, C_\infty(\rit^m,V))$ such that

\medskip\noindent
(i) $\L_{h^\ast}^\rho \circ \chi^u =0$, $i_{h^\ast} \circ \chi^u=0$, $\forall h \in
\rit^m$,

\medskip\noindent
(ii) $(\chi^u(\gamma))(A_0^\ast, \ldots, A_{t-1}^\ast, \alpha_1^\ast, \ldots,
\alpha_u^\ast) = (\gamma(A_0, \ldots, A_{t-1}))(\alpha_1, \ldots, \alpha_u)$
for all $A_i \in g\ell(m,\rit)$ and all $\alpha_j \in \rit^{m\ast}$,

\medskip\noindent
(iii) $(\chi^u(\gamma))(X_0, \ldots, X_{t+u-1})$ is an homogeneous polynomial of degree
\linebreak
$-u+ \sum_i$ degree $(X_i)$, for all $X_i \in s\ell_{m+1}$.
\end{prop}

{\it Proof.} It is easy to verify that the mapping that sends $\gamma$ onto the cochain
\begin{eqnarray*}
\lefteqn{(X_0, \ldots, X_{t+u-1}) \rightarrow} \\
&& \frac{(-1)^t}{t!u!(m+1)^u} \sum_\nu sign(\nu)(\gamma(DX_{\nu_0}, \ldots,
DX_{\nu_{t-1}}))(d\,tr(DX_{\nu_t}), \ldots, d\, tr(DX_{\nu_{t+u-1}}))
\end{eqnarray*}
has the required properties ($\nu$ describes the permutations of $t+u$ elements and
$sign(\nu)$ is the signature of the permutation $\nu$). Hence the existence.

Let $c \in \Lambda^{t+u}(s\ell_{m+1},C_\infty(\rit^m,V))$ be such that $\L_{h^\ast}c=0$,
$i_{h^\ast}c=0$ for each $h \in \rit^m$. Suppose also that $c(X_0, \ldots, X_{t+u-1})$ is
an homogeneous polynomial of degree $-u + \sum_i$ degree $(X_i)$ for all $X_i$. We claim
that if, in addition, 
$$c(A_0^\ast, \ldots, A_{t-1}^\ast, \alpha_1^\ast, \ldots,
\alpha_u^\ast)=0$$
 for all $A_i \in g\ell(m,\rit)$ and all $\alpha_j \in \rit^{m\ast}$,
then $c=0$, thus proving the uniqueness of $\chi^u$. We show that $c(X_0, \ldots,
X_{t+u-1})=0$ by induction on $k=-u + \sum_i$ degree $(X_i)$. It is true for $k \leq 0$,
by assumption. If it is true for $k<\ell$ then for $k=\ell$, it follows from $\L^\rho c=0$
that
$$h^\ast.(c(X_0, \ldots, X_{t+u-1}))=0, \quad \forall h \in \rit^m.$$
The homogeneous polynomial $c(X_0, \ldots, X_{t+u-1})$ of degree $\ell >0$ is thus
constant. Hence it vanishes. \cqfd

\bigskip
Let $\chi$ denote the mapping from $\Lambda(g\ell(m,\rit),\Lambda(\rit^{m\ast},V))$ into
$\Lambda(s\ell_{m+1},C_\infty(\rit^m,V))$ that reduces to $\chi^u$ when restricted to
$\Lambda(g\ell(m,\rit), \Lambda^u(\rit^m,V))$. Observe also that $\Lambda(\rit^{m\ast},V))$ 
is a representation of $g\ell(m,\rit)$, its action being deduced in a natural way from the
given $\rho$ and the natural action on $\rit^{m\ast}$. We denote also by $\partial^\rho$
the Chevalley-Eilenberg coboundary operator associated to that representation.

\begin{lema} \label{L2.1}
Th mapping $\chi$ is an homomorphism of differential spaces.
\end{lema}

{\it Proof.} Let $\gamma$ be given. Then $\L_{h^\ast}^\rho \partial^\rho \chi^u(\gamma) =
\partial^\rho \L_{h^\ast}^\rho \chi^u(\gamma)=0$ and $i_{h^\ast} \partial^\rho
\chi^u(\gamma)=\L_{h^\ast}^\rho \chi^u(\gamma)-\partial^\rho i_{h^\ast}
\chi^u(\gamma)=0$. In addition, it is clear that
$$(\partial^\rho \chi^u(\gamma))(A_0^\ast, \ldots, A_t^\ast, \alpha_1^\ast, \ldots,
\alpha_u^\ast)=((\partial^\rho \gamma)(A_0, \ldots, A_t))(\alpha_1, \ldots, \alpha_u)$$
for all $A_i \in g\ell(m,\rit)$ and all $\alpha_j \in \rit^{m\ast}$. Finally, a careful
examination of the various terms of $(\partial^\rho \chi^u(\gamma))(X_0, \ldots, X_{t+u})$
shows that it is an homogeneous polynomial of degree $-u+\sum_i$ degree $(X_i)$. In view
of the uniqueness property of $\chi^u$, it follows that $\partial^\rho
\chi^u(\gamma)=\chi^u(\partial^\rho \gamma)$. \cqfd

\bigskip
We are now in position to compute the cohomology of $s\ell_{m+1}$ associated to the
representation $L^\rho$.

\begin{theo} \label{T2.3}
The mapping
$$\chi_{\sharp} : H(g\ell(m,\rit), \Lambda(\rit^{m\ast},V)) \rightarrow
H(s\ell_{m+1},C_\infty(\rit^m,V))$$ 
is a bijection.
\end{theo}

{\it Proof.} We use the Hochschild-Serre spectral sequence $E_i^{p,q}$ associated to the
subalgebra $\rit^m$ of $s\ell_{m+1}$. Recall \cite{7} that $E_1^{p,q}$ is the $q$-th
cohomology space of the Chevalley-Eilenberg complex of the representation of $\rit^m$
induced by $h \rightarrow \L_{h^\ast}^\rho$ on the space
$\Lambda^p(s\ell_{m+1}/\rit^m,C_\infty(\rit^m,V))$ or, equivalently, on the space of the
elements of $\Lambda^p(s\ell_{m+1}, C_\infty(\rit^m,V))$ that are cancelled out by each
$i_{h^\ast}$, $h \in \rit^m$.

Once evaluated on ${\bf X} = (X_1, \ldots, X_p) \in s\ell_{m+1}^p$, a $q$-cochain $c$
defines a differential $q$-form on $\rit^m$, namely
$$\omega_{\bf X}^c : (h_1, \ldots, h_q) \rightarrow (c(h_1^\ast, \ldots, h_q^\ast))(X_1,
\ldots, X_p).$$
As easily seen, if $\omega_{\bf X}^c =0$ for degree $({\bf X}) = \sum_i$ degree $(X_i) <
u$, then for \linebreak degree $({\bf X}) =u$, one has
$$d\omega_{{\bf X}}^c = \omega_{{\bf X}}^{\partial^\rho c},$$
where $d$ is the de Rham differential. An easy induction on degree $(X)$ allows thus to
show that

\medskip
(a) {\it If $q>0$, then $E_1^{p,q}=0$.}

\medskip
Now, a $0$-cocycle is an element $c \in \Lambda^p(s\ell_{m+1},C_\infty(\rit^m,V))$ such
that $\L_{h^\ast}^\rho c=0$ and $i_{h^\ast}c=0$ for each $h \in \rit^m$. In particular, if
$c({\bf X})=0$ when degree $({\bf X})<u$, then $c({\bf X})$ is constant for degree $({\bf
X})=u$. Thus, if $\gamma \in \Lambda^{p-u}(g\ell(m,\rit), \Lambda^u(\rit^{m\ast},V))$ is
defined by
$$(\gamma(A_1, \ldots, A_{p-u}))(\alpha_1, \ldots, \alpha_u)=c(A_1^\ast, \ldots,
A_{p-u}^\ast, \alpha_1^\ast, \ldots, \alpha_u^\ast)$$
then $c-\chi^u(\gamma)$ vanishes on arguments ${\bf X}$ of degree $\leq u$. It follows by
induction on $u$ that $c$ belongs to the image of $\chi$. It is clear that $\chi$ is
injective. Therefore, 

\medskip
(b) $\chi : \Lambda(g\ell(m,\rit),\Lambda(\rit^{m\ast},V)) \rightarrow \oplus_{p \geq 0}
E_1^{p,0}$  {\it is an isomorphism of differential spaces.}

\medskip 
We deduce from (a) that the higher differentials $d_i$, $i >1$, of the spectral sequence
are vanishing. The result then follows from (b). \cqfd

\section{Deforming the representations of $g\ell(m,\rit)$}

\indent\indent Let again $(V,\rho)$ be a finite dimensional representation of
$g\ell(m,\rit)$. For each $\lambda \in \rit$, define $\rho_\lambda$ by
$$\rho_\lambda : A \rightarrow \rho(A)- \lambda\,tr(A) id$$
where $tr$ denotes the trace and $id$ the identity from $V$ into $V$.

We want to compute the cohomology of the representation $(V, \rho_\lambda)$ of
$g\ell(m,\rit)$ in terms of that of $(V,\rho)$.

Let $\Phi : (\A,\partial) \rightarrow (\B, \partial)$ be an homomorphism of differential
spaces. The mapping
$$\partial_{\Phi} : (a,b) \mapsto (\partial a, \Phi(a)-\partial b)$$
is a differential on $\A \oplus \B$.

\begin{lema} \label{L3.1}
The cohomology space $H(\A \oplus \B, \partial_\Phi)$ is isomorphic to $\ker \Phi_\sharp
\oplus H(\B,\partial)/\im \,\Phi_\sharp$.
\end{lema}

{\it Proof.} Indeed, the connecting homomorphism (\cite{8}) of the short exact sequence of
differential spaces
$$0 \rightarrow \B \stackrel{i}{\rightarrow} \A \oplus \B \stackrel{j}{\rightarrow} \A
\rightarrow 0$$
where $iy=(0,-y)$ and $j(x,y)=x$, is $\Phi_\sharp$. \cqfd

\begin{prop} \label{P3.2}
The Chevalley-Eilenberg cohomology of $(V, \rho_\lambda)$ is isomorphic to
$$\ker \rho_\lambda({\bf 1})_\sharp \oplus H(s\ell(m,\rit),V)/\im \rho_\lambda(\bf
1)_\sharp$$
where $\rho_\lambda(\bf 1)_\sharp$ is the map induced in cohomology by
$$c \in \Lambda(s\ell(m,\rit),V) \rightarrow \rho_\lambda({\bf 1}) \circ c \in
\Lambda(s\ell(m,\rit),V).$$
\end{prop}

{\it Proof.} Indeed, the map
\begin{equation} \label{E3.1}
c \rightarrow (c|_{s\ell(m,\rit)}, (i_{\bf 1}c)|_{s\ell(m,\rit)})
\end{equation}
is an isomorphism of differential spaces between $\Lambda(g\ell(m,\rit),V)$ equipped with
the differential $\partial^{\rho_\lambda}$ and $(\Lambda s\ell(m,\rit))^2$ equipped with
the differential
$$(c',c'') \rightarrow (\partial^\rho c', \rho_\lambda({\bf 1}) \circ c'-\partial^\rho
c'').$$
Hence the result, in view of Lemma \ref{L3.1}. \cqfd

\bigskip
We shall now apply the above proposition when $V \subset \otimes_b^a \rit^m$ and when
$\rho$ is the natural representation of $g\ell(m,\rit)$ on the tensors. For the sake of
simplicity, $V_\lambda$ will denote the representation $(V,\rho_\lambda)$.

\begin{prop} \label{P3.3}
Let $V$ be a subspace of $\otimes_b^a \rit^m$ stable under the natural representation of
$g\ell(m,\rit)$ on the tensors. If $\lambda \neq (a-b)/m$, then
$H(g\ell(m,\rit),V_\lambda)=0$. If $\lambda = (a-b)/m$, then $H(g\ell(m,\rit),V_\lambda)$
is isomorphic to
\begin{equation} \label{E3.2}
(\Lambda g\ell(m,\rit)^\ast)_{g-inv} \otimes V_{s-inv}
\end{equation}
where $g-inv$ and $s-inv$ denote the invariance with respect to $g\ell(m,\rit)$ and
$s\ell(m,\rit)$ respectively.
\end{prop}

{\it Proof.} If $m>1$, the space of cocycles of the Chevalley-Eilenberg cohomology of the
representation $(V,\rho)$ of $s\ell(m,\rit)$ is the direct sum of the space of
coboundaries and of $(\Lambda s\ell(m,\rit)^\ast)_{s-inv} \otimes V_{s-inv}$ \cite{7}. On
the latter, $\rho_\lambda(\bf 1)_\sharp$ is the multiplication by $a-b-\lambda m$.
Moreover, (\ref{E3.1}) pulls $(\Lambda s\ell(m,\rit)^\ast)_{s-inv}^2$ back onto $(\Lambda
g\ell(m,\rit)^\ast)_{g-inv}$. Hence the result when $m>1$. For $m=1$, the result follows
from an immediate direct computation, (\ref{E3.2}) being just $V$ in this case. \cqfd

\begin{coro} \label{C3.4}
Under the asumptions of Propositon \ref{P3.3}, $H(g\ell(m,\rit),V_\lambda) \neq 0$ only if
$\lambda=(a-b)/m$ is an integer.
\end{coro}

{\it Proof.} It is obvious if $m=1$. If $m>1$, it follows from the description of
$(\otimes_b^a \rit^m)_{s-inv}$ \cite{9}~: $T \in \otimes_b^a \rit^m$ is
$s\ell(m,\rit)$-invariant if and only if 
\linebreak $T(\xi_1, \ldots, \xi_a,X_1, \ldots, X_b)$, as a
function of $\xi_i \in \rit^{m\ast}$ and $X_j \in \rit^m$, is a linear combination of
products of contractions $\scal{X_j}{\xi_i}$ and of determinants of the form
$\det(\xi_{i_1}, \ldots, \xi_{i_m})$ and $\det(X_{j_1}, \ldots, X_{j_m})$. \cqfd

\section{Two examples}

\indent\indent
In this section, we illustrate the previous results by computing two cohomologies of
$s\ell_{m+1}$ that will be needed later.

Let us first consider $S_\lambda^k \rit^m$, the space $\vee^k \rit^m$ of $k$-contravariant
symmetric tensors equipped with the representation $\rho_\lambda$ (where $\rho$ is the
natural representation of $g\ell(m,\rit)$ on tensors).

\begin{prop} \label{P4.1}
(i) Assume that $m>1$. If $(\lambda,k) \in \{(0,0),(1,0),(1,1)\}$, then $H(s\ell_{m+1},
C_\infty (\rit^m,S_\lambda^k \rit^m))$ is isomorphic to $(\Lambda
g\ell(m,\rit)^\ast)_{g-inv}$ otherwise it is vanishing.

\medskip \noindent 
(ii) If $(\lambda,u) \in \{(k,0),(k,1),(k+1,1),(k+1,2)\}$ then
$H(s\ell_2,C_\infty(\rit,S_\lambda^k \rit))$ is isomorphic to $\rit$ otherwise it vanishes.
\end{prop}

{\it Proof.} It follows from Theorem \ref{T2.3} and Proposition \ref{P3.3} that \linebreak 
$H^u(s\ell_{m+1},C_\infty(\rit^m,S_\lambda^k \rit^m))$ is non vanishing only if $\lambda =
(j+k)/m$ for some $j \in \{0, \ldots,m\}$ and that it is then isomorphic to
\begin{equation} \label{E3}
(\Lambda^{u-j} g\ell(m,\rit)^\ast)_{g-inv} \otimes
\Lambda^j(\rit^{m\ast},\vee^k\rit^m)_{s-inv}.
\end{equation}
Using the description of $(\otimes_b^a \rit^m)_{s-inv}$ recalled in the proof of Corollary
\ref{C3.4}, the latter is non vanishing if and only if $(j,k)$ equals (0,0), $(m,0)$ or
$(m-1,1)$ in which case, (\ref{E3}) is respectively spanned by the mappings
\begin{eqnarray*}
(A_0, \ldots, A_{t-1}) & \rightarrow & \gamma(A_0, \ldots, A_{t-1}), \\
(A_0, \ldots, A_{t-1}, \alpha_1, \ldots, \alpha_m) & \rightarrow & \gamma(A_0, \ldots,
A_{t-1}) \det (\alpha_1, \ldots, \alpha_m)
\end{eqnarray*}
or
$$(A_0, \ldots, A_{t-1}, \alpha_1, \ldots, \alpha_{m-1}) \rightarrow [\xi \rightarrow
\gamma(A_0, \ldots, A_{t-1}) \det (\alpha_1, \ldots, \alpha_{m-1},\xi)]$$
$(t=u-j$, $A_i \in g\ell(m,\rit)$, $\alpha_i, \xi \in \rit^{m\ast}$, $\gamma \in
(\Lambda^t g\ell(m,\rit)^\ast)_{g-inv}$). Hence (i).

For (ii), it suffices to note that $\Lambda^j(\rit,\vee^t\rit)_{s-inv}=\rit$. \cqfd

\begin{coro} \label{C4.2}
The cohomology of $s\ell_{m+1}$ acting on the space of scalar $\lambda$-densities of
$\rit^m$ vanishes if $\lambda \neq 0$ and $\lambda \neq 1$ otherwise, it is isomorphic to
$(\Lambda g\ell(m,\rit)^\ast)_{g-inv}$.
\end{coro}

{\it Proof.} It is the cohomology of the representation $C_\infty(\rit^m,S_\lambda^0
\rit^m)$ of $s\ell_{m+1}$. \cqfd

\bigskip
Let us now consider the case of the space $S_p^{k,q}\rit^m$ of $k$-contravariant symmetric
tensors valued in the space of linear maps from $\vee^p\rit^m$ into $\vee^q\rit^m$~:
\begin{eqnarray*}
S_p^{k,q}\rit^m &=& \vee^k \rit^m \otimes Hom(\vee^p \rit^m,\vee^q\rit^m) \\
& \simeq & \vee^k \rit^m \otimes \vee^q \rit^m \otimes \vee^p \rit^{m\ast}.
\end{eqnarray*}

We equip it with the natural representation of $g\ell(m,\rit)$.

\begin{prop} \label{P4.3}
If $p \in \{k+q,k+q+1\}$ then $H(s\ell_{m+1},C_\infty (\rit^m,S_p^{k,q}\rit^m))$ is
isomorphic to $(\Lambda g\ell(m,\rit)^\ast)_{g-inv}$ otherwise it vanishes.
\end{prop}

{\it Proof.} The proof is similar to that of the previous proposition. One is now left to
study
\begin{equation} \label{E4}
(\Lambda^{u-j} g\ell(m,\rit)^\ast)_{g-inv} \otimes
\Lambda^j(\rit^{m\ast},S_p^{k,q}\rit^m)_{s-inv}
\end{equation}
knowing that $p=j+k+q$ (since $``\lambda"=0$). This space is non vanishing if and only if
$j=0$ or $j=1$, (\ref{E4}) being then respectively spanned by the mappings
$$(A_0, \ldots, A_{t-1}) \rightarrow [P \rightarrow \gamma(A_0,\ldots,A_{t-1})(\eta
D_\xi)^kP]$$
or
$$(A_0, \ldots, A_{t-1},\alpha) \rightarrow [P \rightarrow
\gamma(A_0,\ldots,A_{t-1})(\alpha D_\xi)(\eta D_\xi)^kP]$$
($t=u-j$, $A_i \in g\ell(m,\rit)$, $\alpha \in \rit^m$, $\gamma \in (\Lambda^t
g\ell(m,\rit)^\ast)_{g-inv}$, $P \in \vee^k\rit^m$).

Here, we wiew an element of $S_p^{k,q}\rit^m$ as being a homogeneous polynomial of degree
$k$ in $\eta \in \rit^{m\ast}$ valued in the space of linear mappings from the space of
homogeneous polynomials of degree $p$ in $\xi \in \rit^{m\ast}$ into the space of these of
degree $q$. Moreover, $\alpha D_\xi$ and $\eta D_\xi$ denote the derivations with respect
to $\xi$ in the direction $\alpha \in \rit^{m\ast}$ and $\eta \in \rit^{m\ast}$
respectively. These conventions will be helpfull in the sequel. \cqfd

\section{Cohomology of $s\ell_{m+1}$ valued in
$\D(\S_\delta^p(\rit^m),\S_\delta^q(\rit^m))$}

\indent\indent
Recall from the introduction that the structure of the $s\ell_{m+1}$-module
$\D_{\lambda\mu}(\rit^m)$ is related to the cohomology of $s\ell_{m+1}$ with coefficients
in $\D(\S_\delta^p(\rit^m),\S_\delta^q(\rit^m))$, \linebreak $\delta = \mu-\lambda$.

We first compute the cohomology of $s\ell_{m+1}$ with coefficients in
\linebreak $\D^k(\S_\delta^p(\rit^m),\S_\delta^q(\rit^m))$, $k \in \nit$. We proceed by
induction on $k$, using the short exact sequence
\begin{eqnarray*}
0Ê\rightarrow \D^{k-1}(\S_\delta^p(\rit^m),\S_\delta^q(\rit^m)) &
\stackrel{i}{\rightarrow} & \D^k(\S_\delta^p(\rit^m),\S_\delta^q(\rit^m)) \\
&\stackrel{\sigma}{\rightarrow} & C_\infty (\rit^m, S_p^{k,q}\rit^m) \rightarrow 0
\end{eqnarray*}
that induces an exact triangle
\begin{equation} \label{E5}
\begin{array}{cl}
H(s\ell_{m+1},\D^{k-1}(\S_\delta^p(\rit^m),\S_\delta^q(\rit^m))) \\
& \searrow i_\sharp \\
\theta\;\big\uparrow & H(s\ell_{m+1},\D^{k}(\S_\delta^p(\rit^m),\S_\delta^q(\rit^m))) \\
& \swarrow \sigma_\sharp \\
H(s\ell_{m+1},C_\infty(\rit^m,S_p^{k,q}\rit^m))
\end{array}
\end{equation}
where the connecting homomorphism $\theta$ is of degree 1 (\cite{8}). Note that the
cohomology of the quotient is that computed in the previous section. From Proposition
\ref{P4.3}, we immediately deduce

\begin{theo} \label{T5.1}
(a) If $p<q$, then $H(s\ell_{m+1}, \D^{k}(\S_\delta^p(\rit^m),\S_\delta^q(\rit^m)))=0$

\noindent (b) If $n \geq 0$, then

\medskip\noindent (i) $k<n-1 \Rightarrow H(s\ell_{m+1},
\D^{k}(\S_\delta^{q+n}(\rit^m),\S_\delta^q(\rit^m)))=0$

\medskip\noindent (ii) $H(s\ell_{m+1},
\D^{n-1}(\S_\delta^{q+n}(\rit^m),\S_\delta^q(\rit^m)))$ is isomorphic to $(\Lambda
g\ell(m,\rit)^\ast)_{g-inv}$

\medskip\noindent (iii) The inclusion of $
\D^{n}(\S_\delta^{q+n}(\rit^m),\S_\delta^q(\rit^m)))$ into 
$ \D^{k}(\S_\delta^{q+n}(\rit^m),\S_\delta^q(\rit^m))$ induces an isomorphism
$$H(s\ell_{m+1}, \D^{n}(\S_\delta^{q+n}(\rit^m),\S_\delta^q(\rit^m))) \rightarrow
H(s\ell_{m+1} , \D^{k}(\S_\delta^{q+n}(\rit^m),\S_\delta^q(\rit^m))) $$
for each $k>n$.
\end{theo}

\begin{rema} \label{R5.2}{\rm In (ii), the isomorphism is obtained by composing the
isomorphism of Propositon \ref{P4.3} and the map induced in cohomology by the symbol map
$\sigma$ from
$\D^{n-1}(\S_\delta^{q+n}(\rit^m), \S_\delta^q(\rit^m))$.}
\end{rema}

To complete Theorem \ref{T5.1}, we need to compute the connecting homomorphism~$\theta$.

\begin{lema} \label{L5.2}
If $n>0$, identifying $H(s\ell_{m+1},C_\infty(\rit^m,S_{q+n}^{n,q},\rit^m))$ and
\linebreak 
$H(s\ell_{m+1}, \D^{n-1}(\S_\delta^{q+n}(\rit^m),\S_\delta^q(\rit^m))) $ with $(\Lambda
g\ell(m,\rit)^\ast)_{g-inv}$ using Proposition \ref{P4.3} and Theorem \ref{T5.1}, one has
$$\theta(\gamma) = (-1)^{t+1}n[(m+1)\delta - (2q+n+m)]\gamma,\quad \forall \gamma \in
(\Lambda^t g\ell(m,\rit)^\ast)_{g-inv}.$$
\end{lema}

{\it Proof.} To perform the computation, it is convenient to identify the spaces
$\D^k(\S_\delta^p(\rit^m),\S_\delta^q(\rit^m))$ and $\oplus_{0 \leq i \leq k}
C_\infty(\rit^m, S_p^{i,q} \rit^m)$ by representing the operator
$$D : T \rightarrow \sum_{\mod{\alpha} \leq k} A_\alpha(D_x^\alpha T)$$
where $A_\alpha \in C_\infty(\rit^m,Hom(\vee^p\rit^m,\vee^q\rit^m))$, by
$$\T_D : (\eta,P) \in \rit^{m\ast} \times \vee^p \rit^m \rightarrow \sum_{\mod{\alpha} \leq
k} \eta^\alpha A_\alpha(P) \in C_\infty(\rit^m, \vee^q \rit^m).$$
(cf. the proof of Proposition \ref{P4.3}). It follows then easily that, the components of
$\zeta \in \rit^{m\ast}$ representing symbolically the partial derivatives of $X \in
Vect(\rit^m)$,
\begin{eqnarray}
\lefteqn{\T_{L_XD}(\eta,P)} \nonumber \\ &=& (X.\T_D)(\eta,P) + \scal{X}{\eta}\T_D(\eta,P)
+\delta\scal{X}{\zeta}\T_D(\eta,P) \nonumber\\
&&-X(\zeta D_\xi)\T_D(\eta,P)-
\T_D(\eta+\zeta,\scal{X}{\eta}P+\delta\scal{X}{\zeta}P-X(\zeta D_\xi)P) \label{E6}
\end{eqnarray}
(see \cite{5} where similar symbolical computations are used).

Now $\gamma \in (\Lambda g\ell(m,\rit)^\ast)_{g-inv}$ defines an element in
$H(s\ell_{m+1},C_\infty(\rit^m,S_{q+n}^{n,q} \rit^m))$, namely the class of the cocycle
$c=\chi(\gamma)$
$$(X_0, \ldots, X_{t-1}) \mapsto (-1)^t \gamma(DX_0, \ldots, DX_{t-1})(\eta D_\xi)^n.$$

In view of the above identification, $c$ may be considered as valued in \linebreak 
$\D^n(\S_\delta^{q+n}(\rit^m)$, $\S_\delta^q(\rit^m))$. Its coboundary is valued in 
$\D^{n-1}(\S_\delta^{q+n}(\rit^m),\S_\delta^q(\rit^m))$ and we seek for the $\gamma' \in
(\Lambda g\ell(m,\rit)^\ast)_{g-inv}$ which $[\partial c]$ corresponds to through the
isomorphism (ii) of Theorem \ref{T5.1} (see Remark \ref{R5.2}). To compute $\gamma'$, it
suffices to evaluate $(\partial c)(A_0^\ast, \ldots, A_{t-1}^\ast, \alpha^\ast)$, $A_i \in
g\ell(m,\rit)$, $\alpha \in \rit^m$, at $x=0$. That is the same to compute the constant
term of
\begin{equation} \label{E7}
(-1)^t L_{\alpha^\ast} (c(A_0^\ast, \ldots, A_{t-1}^\ast))
\end{equation}
where $ L_{\alpha^\ast}$ is the derivation of differential operators in the direction of
$\alpha^\ast$. Observe that
$$c(A_0^\ast, \ldots, A_{t-1}^\ast) = \gamma (A_0, \ldots, A_{t-1})(\eta D_\xi)^n$$
has constant coefficients so that we need only the terms of (\ref{E7}) involving second
order derivatives of the coefficients of $\alpha^\ast$. It follows from (\ref{E6}) that if
$X=\alpha^\ast$, the terms of second order in $X$ in $\T_{L_XD}(\eta,P)$ are given by
\begin{eqnarray*}
\lefteqn{-(\eta D_\eta)(\alpha D_\eta) \T_D(\eta,P) -(m+1)\delta(\alpha
D_\eta)\T_D(\eta,P)}\\
&&+(\alpha D_\eta)\T_D(\eta,(\xi D_\xi)P)+ \sum_i
D_{\eta_i}\T_D(\eta,\xi_i(\alpha D_\xi)P).
\end{eqnarray*}

With $\T_D(\eta,P)=(\eta D_\xi)^nP$ and $P \in \vee^{q+n}\rit^m$, this gives 
$$-n[(m+1)\delta-(2q+n+m)](\alpha D_\xi)(\eta D_\xi)^{n-1}.$$
Hence the lemma. \cqfd

\begin{theo} \label{T5.3}
If $n>0$, the space $H(s\ell_{m+1}, \D^n(\S_\delta^{q+n} (\rit^m),\S_\delta^q(\rit^m)))$
is isomorphic to $(\Lambda g\ell(m,\rit)^\ast)^2_{g-inv}$ if $\delta=(2q+n+m)/(m+1)$ and
vanishes otherwise. If $n=0$, it is isomorphic to $(\Lambda g\ell(m,\rit)^\ast)_{g-inv}$.
\end{theo}

{\it Proof.} This follows immediately from (\ref{E5}) and previous results. \cqfd

\bigskip
We are now able to compute the cohomology of $s\ell_{m+1}$ with coefficients in \linebreak 
$\D(\S_\delta^p(\rit^m),\S_\delta^q(\rit^m))$.

\begin{theo} \label{T5.4}
If $p<q$, then $H(s\ell_{m+1}, \D(\S_\delta^p(\rit^m),\S_\delta^q(\rit^m)))=0$. Otherwise,
the inclusion
$$\D^n(\S_\delta^p(\rit^m),\S_\delta^q(\rit^m)) \stackrel{i}{\rightarrow}
\D(\S_\delta^p(\rit^m),\S_\delta^q(\rit^m))$$
induces an isomorphism in cohomology, where $n=p-q\geq 0$.
\end{theo}

{\it Proof.} For the sake of simplicity, we denote here by $\D$ and $\D^k$ the spaces
$\D(\S_\delta^p(\rit^m),\S_\delta^q(\rit^m))$ and
$\D^k(\S_\delta^p(\rit^m),\S_\delta^q(\rit^m))$ respectively.

Observe that a cochain on $s\ell_{m+1}$ valued in the space of differential operators
necessarily takes its values in the space of operators of a certain order $k$. This is
because $s\ell_{m+1}$ is finite dimensional.

This implies first that $H(s\ell_{m+1},\D)=0$ if $p<q$ because, in this case, the spaces 
$H(s\ell_{m+1},\D^k)$ are all vanishing.

Second, this also implies that if $p=q+n$, $n \geq 0$,
$$i_\sharp : H(s\ell_{m+1},\D^n) \rightarrow H(s\ell_{m+1},\D)$$
is an isomorphism. Indeed, if $i_\sharp [S]=0$, then $S$ is the coboundary of some
$\D^k$-valued cochain $T$. If $k \leq n$, then $[S]=0$ because $\D^k \subset \D^n$. If
$k>n$, then again $[S]=0$, due to Theorem \ref{T5.1} (b) (iii). This proves that
$i_\sharp$ is injective. A similar discussion shows that it is onto. \cqfd

\begin{rema} \label{R5.5}{\rm
When $p>q$, the cohomology of $s\ell_{m+1}$ with coefficients in \linebreak 
$\D(\S_\delta^p(\rit^m),\S_\delta^q(\rit^m))$ is non vanishing if and only if $\delta =
\frac{m+p+q}{m+1}$. This is why we say that $\frac{m+p+q}{m+1}$ is a {\it critical value}
for $\delta$. When $m=1$, these critical values are the special values pointed out in
\cite{3} in classifying the $s\ell_2$-module $\D_{\lambda \mu}^k(\rit)$. In \cite{3},
they were called ``resonant" instead of critical.}
\end{rema}

For $n \geq 0$, we define
$$\tau_n : s\ell_{m+1} \rightarrow \D(\S_\delta^{\ast +n}(\rit^m),\S_\delta^\ast(\rit^m))$$
by
$$\tau_n (X)P = - \sum_{j,i_1, \ldots, i_n} \partial_jX^j \partial_{i_1} \ldots
\partial_{i_n} D_{\xi_{i_1}} \ldots  D_{\xi_{i_n}}  P$$
(recall that $\partial_i$ denotes the partial derivation with respect to the $i$-th
coordinate and that $D_{\xi_j}$ denotes the derivation with respect to the $j$-th
component of $\xi \in \rit^{m\ast}$). 

Using the notations of the proof of Lemma \ref{L5.2}, one can also write
$$\tau_n(X) =- \scal{X}{\zeta}(\eta D_\xi)^n$$
where $\eta$ and $\zeta$ represent the derivatives acting on $P$ and $X$ respectively.

For $n>0$, we also introduce the map
$$\gamma_n : s\ell_{m+1} \rightarrow \D(\S_\delta^{\ast+n}(\rit^m),
\S_\delta^\ast(\rit^m))$$
given by
$$\gamma_n(X)P = \frac{1}{m+1} \sum_{j,i_1, \ldots, i_n} \partial_{i_1} \partial_j X^j
\partial_{i_2} \ldots \partial_{i_n}  D_{\xi_{i_1}}  \ldots  D_{\xi_{i_n}} P.$$
It is symbolically given by
$$\gamma_n(X) = \frac{1}{m+1} \scal{X}{\zeta} (\zeta D_\xi)(\eta D_\xi)^{n-1}.$$

\begin{coro} \label{C5.6}
Let $p=q+n$. The space $H^1(s\ell_{m+1},\D(\S_\delta^p(\rit^m),\S_\delta^q(\rit^m))$ is
spanned by $[\tau_0]$ if $n=0$ and by $[\tau_n]$ and $[\gamma_n]$ if $n>0$ and $\delta =
\frac{m+p+q}{m+1}$. It vanishes otherwise.
\end{coro}

{\it Proof.} This easily follows from the previous results. \cqfd

\section{The space $Hom_{s\ell_{m+1}} (\S_\delta^p(\rit^m),\S_\delta^q(\rit^m))$}

\indent\indent
In this section, we compute the $0$-th cohomology space of $s\ell_{m+1}$ with
coefficients in $Hom(\S_\delta^p(\rit^m),\S_\delta^q(\rit^m))$. It is the space
$Hom_{s\ell_{m+1}} (\S_\delta^p(\rit^m),\S_\delta^q(\rit^m))$ of $s\ell_{m+1}$-equivariant
linear mappings from $\S_\delta^p(\rit^m)$ into $\S_\delta^q(\rit^m)$.

For each $n \in \nit$, we define
$$T_n : \S_\delta^{\ast+n}(\rit^m) \rightarrow \S_\delta^\ast(\rit^m)$$
by
$$T_nP = \sum_{i_1, \ldots, i_n} \partial_{i_1} \ldots \partial_{i_n}  D_{\xi_{i_1}} 
\ldots  D_{\xi_{i_n}} P.$$
With the notations of the proof of Lemma \ref{L5.2}, $T_n$ is represented by the
polynomial $(\eta D_\xi)^n$. Moreover, it follows from that proof that
$$L_X \circ T_n - T_n \circ L_X = -n [(m+1)\delta -(m+2p+n)] \gamma_n(X).$$
In particular, $T_n : \S_\delta^{p+n}(\rit^m) \rightarrow \S_\delta^p(\rit^m)$ {\it is
$s\ell_{m+1}$-equivariant if and only if} \linebreak $n[(m+1)\delta- (m+2p+n)]=0$.

Observe that
$$\tau_n(X)=-\scal{X}{\zeta} T_n$$
and 
$$\gamma_n(X)= \frac{1}{m+1} \scal{X}{\zeta}(\zeta D_\xi) \circ T_{n-1}.$$
Moreover,
$$\gamma_{n+i}(X)= \gamma_n(X) \circ T_i$$
since $T_i \circ T_j = T_{i+j}$.

\begin{lema} \label{L6.1}
If $p\geq q$ and if $A \in Hom(\S_\delta^p(\rit^m),\S_\delta^q(\rit^m))$ commute with the
Lie derivations in the direction of $\partial_1, \ldots, \partial_m$ and of
$E=x^i\partial_i$ then $A$ is a differential operator.
\end{lema}

{\it Proof.} It has been shown in \cite{4} that if $\delta=0$, then under the assumptions
of the lemma, $A$ is a local operator. In fact, the proof of \cite{4} works also if
$\delta \neq 0$. Thus, we may assume that $A$ is local. From the theorem of Peetre
\cite{10}, it is then locally a differential operator. Since it commutes with
$L_{\partial_1}, \ldots,L_{\partial_m}$, its order is bounded. \cqfd

\begin{theo} \label{T6.2}
If $p \geq q$ and if $A \in Hom_{s\ell_{m+1}} (\S_\delta^p(\rit^m),\S_\delta^q(\rit^m))$
is non vanishing then either $p=q$ and $A$ is a constant multiple of the identity or
$p>q$, $\delta = \frac{m+p+q}{m+1}$ and $A$ is a constant multiple of $T_{p-q}$.
\end{theo}

{\it Proof.} Indeed, since $A$ is a differential operator, it is a $0$-cocycle of
$s\ell_{m+1}$ valued in $\D(\S_\delta^p(\rit^m),\S_\delta^q(\rit^m))$. The theorem then
follows immediately from Theorem~\ref{T5.4}. \cqfd

\bigskip
The case $p<q$ is more difficult. We state directly the theorem. Indeed, we do not know
whether a lemma similar to Lemma \ref{L6.1} holds true in this case.

\begin{theo} \label{T6.3}
If $p<q$ and $A \in Hom_{s\ell_{m+1}} (\S_\delta^p(\rit^m),\S_\delta^q(\rit^m))$, then
$A=0$.
\end{theo}

{\it Proof.} It suffices to show that $A$ is local (and thus a differential operator,
proceeding like in the proof of Lemma \ref{L6.1}). Indeed, it follows from Theorem
\ref{T5.4} that a $0$-cocycle of $s\ell_{m+1}$ with values in
$\D(\S_\delta^p(\rit^m),\S_\delta^q(\rit^m))$ is vanishing, because $p<q$.

The proof has three parts. In the first, $\delta$ has not the cricital value
$\frac{m+p+q}{m+1}$ and we show directly that $A=0$. In the second, $\delta =
\frac{m+p+q}{m+1}$ and $m>1$ while in the third, $\delta = \frac{m+p+q}{m+1}$
but $m=1$. We set again $q=p+n$.

\medskip
(i) {\it The non critical case.} We make use of the Casimir operator $\C_\delta^p$ of the
$s\ell_{m+1}$-module $\S_\delta^p(\rit^m)$. Recall from \cite{11} that it is the
$s\ell_{m+1}$-equivariant linear map from $\S_\delta^p(\rit^m)$ into itself defined by
$$\C_\delta^p = \sum_i L_{X_i^\ast} \circ L_{Y_i^\ast}$$
where $\{X_i :i \leq m(m+2)\}$ is any basis of $s\ell(m+1,\rit)$ and $\{Y_i : i \leq
m(m+2)\}$ is its dual with respect to the Killing form $\bf K$ of $s\ell(m+1,\rit)$ (i.e.
${\bf K}(X_i,Y_j)=\delta_{ij})$. By Theorem \ref{T6.2}, $\C_\delta^p$ is a multiple
$c_\delta^p$ of the identity. It is easy to compute. One gets
$$2(m+1)c_\delta^p = m(m+1)\delta^2 -(m+2p)(m+1)\delta + 2p(m+p).$$
Since $A$ is $s\ell_{m+1}$-equivariant, one has $A \circ \C_\delta^p = \C_\delta^{p+n}
\circ A$. But
$$c_\delta^{p+n} - c_\delta^p =-n\left(\delta-\frac{m+2p+n}{m+1} \right).$$
Therefore, if $\delta \neq \frac{m+p+q}{m+1}$, then $A=0$.

\medskip
(ii) {\it The critical case} $\delta = \frac{m+2p+n}{m+1}$, $m>1$. This case is more
delicate to handle. We need some preparations. For each $\alpha \in \rit^{m\ast}$, $\hat
\alpha$ denotes the function $x \rightarrow \alpha(x)$ on $\rit^m$ (it is $\alpha$ but
viewed as an element of $C_\infty(\rit^m)$). Moreover, $I \in \S_0^1(\rit^m)$ is defined by
$$I(\xi)_x = \xi(x),\quad \forall x \in \rit^m,\;\;\forall \xi \in \rit^{m\ast}.$$

One has
$$L_{\alpha^\ast}Q = L_E(\hat \alpha Q) +(\delta-1) \hat \alpha Q -I(\alpha D_\xi)Q, \quad
\forall Q \in \S_\delta^\ast(\rit^m).$$
Using this, the $s\ell_{m+1}$-equivariance of $A$ and the assumption $\delta =
\frac{m+2p+n}{m+1}$, one gets, after some computations,
$$(m+p-1) A(\hat \alpha Q)_0 = A(I(\alpha D_\xi)Q)_0, \quad \forall Q \in
\S_\delta^p(\rit^m),$$
the index $0$ denoting the evaluation at $x=0$. This implies easily that
$$A(\hat \alpha I^kQ)_0 = \frac{1}{m+p-(k+1)}\;A(I^{k+1}(\alpha D_\xi)Q)_0.$$
Hence, by induction on $k$,
$$A(\hat \alpha_1 \ldots \hat \alpha_kQ)_0 = \frac{1}{(m+p-1) \ldots
(m+p-k)}\;A(I^k(\alpha_1D_\xi) \ldots (\alpha_kD_\xi)Q)_0$$
for each $Q \in \S_\delta^p(\rit^m)$ and for all $\alpha_1, \ldots,\alpha_k \in
\rit^{m\ast}$.

We now prove that $A$ is local. Suppose that $P \in \S_\delta^p(\rit^m)$ vanishes in a
neighborhood of $y \in \rit^m$. We want to show that $A(P)_y=0$.

Assume first that $y=0$. There exists then $Q_{i_0 \ldots i_p} \in \S_\delta^p(\rit^m)$
such that 
$$P = \widehat{\varepsilon^{i_0}} \ldots \widehat{\varepsilon^{i_p}} Q_{i_0 \ldots i_p}$$
($\varepsilon^1, \ldots, \varepsilon^m$ denotes canonical basis of $\rit^{m\ast}$ :
$\varepsilon^i(x) = x^i,$ $\forall x \in \rit^m$). Then
$$A(P)_0 = \frac{1}{(m+p-1) \ldots (m-1)}\;A(I^{p+1}  D_{\xi_{i_0}}  \ldots 
D_{\xi_{i_p}}  Q_{i_0
\ldots i_p})_0=0$$
because the $Q_{i_0 \ldots i_p}$'s are homogeneous of order $p$ in $\xi \in \rit^{m\ast}$.

One reduces the case $y \neq 0$ to the case $y=0$ just by replacing above $x$ by $x-y$ in
$\hat \alpha,I,E,\alpha^\ast$, etc.

\medskip
(iii) {\it The critical case} $\delta = \frac{m+2p+n}{m+1}$, $m=1$. We denote by $t$ the
canonical coordinate of $\rit$. The map $A$ is of the form
$$A : f \left( \frac{d}{dt} \right)^p \mod{dt}^\delta \rightarrow A_f \left( \frac{d}{dt}
\right)^{p+n} \mod{dt}^\delta, \quad f,A_f \in C_\infty(\rit).$$
Expressing the fact that $A$ commutes with $L_{\frac{d}{dt}}$ and $L_{t\frac{d}{dt}}$ leads
immediatelty to the following relations :
$$\frac{d}{dt}\;A_f = A_{\frac{df}{dt}}$$
and
$$t \frac{d}{dt} A_f = A_{t \frac{df}{dt}} + nA_f.$$

From this, it follows first that $A_f=0$ if $f$ is a polynomial. Indeed, the first
relation shows that if $f$ is a homogeneous polynomial of degree $s$, then $A_f$ is a
polynomial of degree $\leq s$ but the second relation shows that, in the same time, $A_f$
is homogeneous of order $s+n >s$.

On the other hand, $T_n : \S_\delta^{p+n}(\rit) \rightarrow \S_\delta^p(\rit)$ is
$s\ell_2$-equivariant (because $\delta$ is critical). By Theorem \ref{T6.2}, it follows
that $A \circ T_n$ is a constant multiple $a$ of the identity on $\S_\delta^{p+n}(\rit)$.
It is clear that $a=0$ because $A$ vanishes on polynomials. This reads
$$A_{\frac{d^n}{dt^n}\, f}=0$$
for all $f \in C_\infty(\rit)$. Therefore, $A=0$. \cqfd

\section{The splitting of the short exact sequence (\ref{E0}) for $M=R^m$}

\indent\indent
Our aim in this section is to obtain a necessary and sufficient condition for the short
exact sequence of $s\ell_{m+1}$-modules (\ref{E0}) to be split (with $M=\rit^m$).

We denote by $\varphi : \S_\delta(\rit^m) \rightarrow \D_{\lambda\mu}(\rit^m)$ the
canonical right inverse of the symbol map. If $P \in \S_\delta^k(\rit^m)$, then
$\varphi(P) \in \D_{\lambda\mu}^k(\rit^m)$ is the unique homogeneous differential operator
of order $k$ such that $\sigma(\varphi(P))=P$.  As known, the coboundary $E_k$ of the
restriction of $\varphi$ to $\S_\delta^k(\rit^m)$ takes its values in
$Hom(\S_\delta^k(\rit^m),\D_{\lambda\mu}^{k-1}(\rit^m))$. Moreover, its cohomology class
characterizes the isomorphism class of the short exact sequence (\ref{E0}) \cite{11}. In
particular, (\ref{E0}) is split if and only if $E_k$ is a coboundary of some element of
$Hom(\S_\delta^k(\rit^m),\D_{\lambda\mu}^{k-1}(\rit^m))$.

\begin{lema} \label{L7.1}
One has $E_k=-u_k \varphi \circ \gamma_1$, where
$$u_k = (m+1) \lambda +k-1.$$
\end{lema}

{\it Proof.} It is just a matter of simple computation. \cqfd

\begin{lema} \label{L7.2}
As a cocycle of $s\ell_{m+1}$ valued in
$Hom(\S_\delta^k(\rit^m),\S_\delta^{k-n}(\rit^m))$, $\gamma_n$ is a coboundary if and only
if
$$v_n =-n((m+1)\delta-(m+2k-n))$$
is non vanishing.
\end{lema}

{\it Proof.} By Corollary \ref{C5.6}, it suffices to show that if $\gamma_n=\partial T$,
then \linebreak $T \in Hom(\S_\delta^k(\rit^m),\S_\delta^{k-n}(\rit^m))$ is a differential
operator. Since $\gamma_n(\partial_i)=0$ and $\gamma_n(E)=0$, this follows immediately
from Lemma
\ref{L6.1}. \cqfd

\begin{rema} \label{R7.3}{\rm
It is clear that if $v_n \neq 0$, then $\gamma_n$ is the coboundary of
$\frac{1}{v_n}\;T_n$.}
\end{rema}

In the sequel, for $\A=\sum_i \A_i t^i \in \rit[t]$, we denote by $\A(\eta\D_\xi)$ the
differential opertor
$$\sum_i \A_i T_i : \S_\delta^k(\rit^m) \rightarrow \bigoplus_{o \leq k}
\S_\delta^{k-i}(\rit^m).$$

\begin{prop} \label{P7.4}
If $\gamma_1, \ldots, \gamma_{n-1}$ are coboundaries, then
\begin{equation} \label{E8}
E_k = - \frac{u_k \ldots u_{k-n+1}}{v_1 \ldots v_{n-1}}\; \varphi \circ \gamma_n +
\partial \varphi \circ \A_n(\eta D_\xi)
\end{equation}
for some polynomial $\A_n$ of degree $n-1$ such that $\A_n(0)=0$.
\end{prop}

{\it Proof.} We proceed by induction on $n$. By Lemma \ref{L7.1}, $E_k = -u_k \varphi \circ
\gamma_1$. Assume that (\ref{E8})  holds true and that $\gamma_n$ is a coboundary. Then,
setting
$$a=-\frac{u_k \ldots u_{k-n+1}}{v_1 \ldots v_{n-1}}$$
for simplicity, one successively gets
\begin{eqnarray*}
\lefteqn{(E_k-\partial \varphi \circ \A_n(\eta D_\xi))(X)} \\
 &=& \frac{a}{v_n}\;(\varphi
\circ L_X
\circ T_n - \varphi \circ T_n \circ L_X) \\
&=& \frac{a}{v_n}\;(\varphi \circ L_X\circ T_n - L_X \circ \varphi \circ T_n) +
\frac{a}{v_n} ( L_X \circ \varphi \circ T_n - \varphi \circ T_n \circ L_X) \\
&=& - \frac{a}{v_n}\; E_{k-n} (X) \circ T_n + \frac{a}{v_n}\;(\partial(\varphi \circ
T_n))(X) \\
&=& \frac{a u_{k-n}}{v_n}\; \varphi \circ \gamma_1(X) \circ T_n +
\frac{a}{v_n}\;(\partial(\varphi \circ T_n))(X) \\
&=& \frac{a u_{k-n}}{v_n}\; \varphi \circ \gamma_{n+1}(X) +
\frac{a}{v_n}\;(\partial(\varphi \circ T_n))(X)
\end{eqnarray*}
because $\gamma_1 \circ T_n = \gamma_{n+1}$. \cqfd

\begin{theo} \label{T7.5}
The short exact sequence of $s\ell_{m+1}$-modules
$$0 \rightarrow \D_{\lambda\mu}^{k-1} (\rit^m) \stackrel{i}{\rightarrow}
\D_{\lambda\mu}^k(\rit^m) \stackrel{\sigma}{\rightarrow} \S_\delta^k \rightarrow 0\quad (k
\geq 1)$$
is split if and only if
\begin{equation} \label{E8.1}
\delta \not\in \left\{ \frac{m+k}{m+1}, \frac{m+k+1}{m+1}, \ldots, \frac{m+2k-1}{m+1}
\right\}
\end{equation} 
or
\begin{equation} \label{E9}
\delta = \frac{m+2k-n}{m+1} \quad and\quad \lambda = \frac{i-k}{m+1}
\end{equation}
for some $n \in \{1, \ldots,k\}$ and some $i \in \{1, \ldots,n\}$.
\end{theo}

{\it Proof.} Let $n$ denote the least integer such that $\gamma_n$ is not a coboundary, as
a cocycle valued in $Hom(\S_\delta^k(\rit^m),\S_\delta^{k-n}(\rit^m))$ (we set $n=k+1$ if
$\gamma_1, \ldots, \gamma_k$ are coboundaries).

It follows from Lemma \ref{L7.2} and Proposition \ref{P7.4} that $E_k$ is a coboundary if
(\ref{E8.1}), or (\ref{E9}), holds true.

Conversely, assume that $E_k$ is a coboundary and that $n \leq k$. Then, $\delta =
\frac{m+2k-n}{m+1}$. Moreover, it follows from (\ref{E7}) that $u_k \ldots u_{k-n+1}
\gamma_n$ is a coboundary. Since $\gamma_n$ is not a coboundary, $u_i=0$ for some $i \in
\{k-n+1, \ldots, k\}$. Hence the result. \cqfd

\begin{coro} \label{C7.6}
If $(m+1)\delta-m \not \in \nit_0$, then there eixsts a unique $s\ell_{m+1}$-equivariant
linear bijection $\sigma^{\lambda\mu}$ : $\D_{\lambda\mu}(\rit^m) \rightarrow
\S_\delta(\rit^m)$ such that, for each $A \in \D_{\lambda\mu}^k(\rit^m)$, the term of
highest order of $\sigma^{\lambda\mu}(A)$ is the symbol $\sigma(A)$ of $A$.
\end{coro}

{\it Proof.} The existence of $\sigma^{\lambda\mu}$ follows immediately from Theorem
\ref{T7.5}. For the uniqueness, assume that $\sigma'$ has the same properties than
$\sigma^{\lambda\mu}$. Then, for each $k$, the restriction of $S=\sigma' \circ
(\sigma^{\lambda\mu})^{-1}$ to $\S_\delta^k(\rit^m)$ is of the form
$$P \in \S_\delta^k(\rit^m) \rightarrow (P,S_1(P), \ldots, S_k(P)) \in \bigoplus_{0 \leq
\ell \leq k} \S_\delta^\ell(\rit^m)$$
for some $s\ell_{m+1}$-equivariant $S_i \in
Hom(\S_\delta^k(\rit^m),\S_\delta^{k-i}(\rit^m))$, $i=1, \ldots,k$. By Theorem \ref{T6.2},
$S_i=0$ for $i=1, \ldots,k$. \cqfd

\begin{rema} \label{R7.7}{\rm
The existence and uniqueness of $\sigma^{\lambda\mu}$ was shown in \cite{4} for
$\lambda=\mu$ and $m>1$. It has been obtained in \cite{3,6} for $m=1$ and non critical
values of $\delta$. For critical $\delta$, $\D_{\lambda\mu}(\rit)$ has been described in
\cite{3}. We will not discuss here the case of the critical values of $\delta$ in higher
dimension. The next result is also a generalization of a useful result of \cite{3,4}.
}\end{rema}

\begin{prop} \label{P7.8}
If $(m+1)\delta-m\not\in \nit_0$, then $T : \D_{\lambda\mu}(\rit^m) \rightarrow
\D_{\lambda'\mu'}(\rit^m) $ is $s\ell_{m+1}$-equivariant if and only if there exists
constant $a_k \in \rit$ such that, for each $k \in \nit$,
$$\sigma^{\lambda'\mu'} \circ T \circ (\sigma^{\lambda\mu})^{-1} = a_k\,id$$
on $\S_\delta^k(\rit^m)$.
\end{prop}

{\it Proof.} This follows immediately from Theorem \ref{T6.2} and Theorem \ref{T6.3}. \cqfd

\section{Cohomology of $s\ell_{m+1}$ valued in $\D_{\lambda\mu}(\rit^m)$}

\indent\indent
The computation of $H(s\ell_{m+1},\D_{\lambda\mu}(\rit^m))$ is quite similar to that of
\linebreak 
$H(s\ell_{m+1},\D(\S_\delta^p(\rit^m),\S_\delta^q(\rit^m)))$. It uses the short exact
sequence (\ref{E0}) and Proposition \ref{P4.1} to get first the spaces 
 $H(s\ell_{m+1},\D_{\lambda\mu}^k(\rit^m))$, $k \in \nit$. In most of the cases, the
cohomology of the module $C_\infty(\rit^m,S_\delta^k\rit^m)$ vanishes. It is however
necessary to compute the connecting homomorphism associated to the sequence (\ref{E0}) for
some values of $\delta$ and $k$. This leads also to a sort of critical values for $\lambda$
and $\mu$. In view of Proposition \ref{P4.1}, they are not the same when $m>1$ as when
$m=1$.

We only summarize the results, leaving the reader to supply the proofs.

\bigskip
A. {\it The case $m>1$}

\begin{theo} \label{T8.1} $\{m>1\}$ If $\delta \not\in \{0,1\}$ or if $\delta=1$ and
$\lambda \neq 0$, then \linebreak 
$H(s\ell_{m+1},\D_{\lambda\mu}(\rit^m))=0$.
\end{theo}

\begin{theo} \label{T8.2} $\{m>1\}$ The space
$H(s\ell_{m+1},\D_{\lambda\lambda}^0(\rit^m))$ is isomorphic to \linebreak 
$(\Lambda
g\ell(m,\rit)^\ast)_{g-inv}$. The inclusion of $\D_{\lambda\lambda}^0(\rit^m)$ into
$\D_{\lambda\lambda}(\rit^m)$ induces an isomorphism in cohomology.
\end{theo}

\begin{theo} \label{T8.3}
$\{m >1\}$ One has a short exact sequence
$$0 \rightarrow H(s\ell_{m+1},\D_{01}^0(\rit^m))\stackrel{i_\sharp}{\rightarrow}
H(s\ell_{m+1},\D_{01}^1(\rit^m))\stackrel{\sigma_\sharp}{\rightarrow} H(s\ell_{m+1},
C_\infty(\rit^m,S_1^1\rit^m)) \rightarrow 0$$
where the kernel and the quotient are isomorphic to $(\Lambda g\ell(m,\rit)^\ast)_{g-inv}$.
The inclusion of $\D_{01}^1(\rit^m)$ into $\D_{01}(\rit^m)$ induces an isomorphism in
cohomology.
\end{theo}

\begin{rema} \label{R8.4}{\rm
If $\lambda \neq 0$, then $H(s\ell_{m+1},\D_{\lambda,\lambda+1}^0(\rit^m))$ is isomorphic
to \linebreak
$(\Lambda g\ell(m,\rit)^\ast)_{g-inv}$ while 
$H(s\ell_{m+1},\D_{\lambda,\lambda+1}^1(\rit^m))=0.$}
\end{rema}

\bigskip
B. {\it The case} $m=1$

\medskip
In this case, the ``critical" values are given by 
$(\lambda,\mu)=\left( \frac{1-n}{2},\frac{1+n}{2} \right)$,  $n \in \nit_0$. Moreover, the
description of the cohomology is very close to that of the module
$\D(\S_\delta^p(\rit),\S_\delta^q(\rit))$. This follows from the fact that this module is
isomorphic to $\D_{\delta-p,\delta-q}(\rit)$ since $\S_\delta^p(\rit)$ is isomorphic to the
space of $(\delta-p)$-densities. In particular, the critical values introduced in Section 6
occur when $p-q$ is a positive integer $n$. They are then given by $\delta=(p+q+1)/2$. In
this case,
$$\delta-p=\frac{1-n}{2}\quad \mbox{and} \quad \delta-q=\frac{1+n}{2}.$$
Note that these critical values were also obtained in \cite{3}, but in a completely
different way.

\begin{theo}\label{T8.51}
If $\delta \not\in \nit$ or if $\delta \in \nit_0$ and $\lambda \neq \frac{1-\delta}{2}$,
then $H(s\ell_2,\D_{\lambda \mu}(\rit))=0$.
\end{theo}

\begin{theo} \label{T8.5}
The inclusion of $\D_{\lambda\lambda}^0(\rit)$ into $\D_{\lambda\lambda}(\rit)$ induces an
isomorphism in cohomology. In particular
$$H^i(s\ell_2, \D_{\lambda\lambda}(\rit)) = \left\{ \begin{array}{lcl}
\rit & \mbox{if} & i=0,1 \\
0 & \mbox{if} & i=2,3
\end{array} \right.$$
\end{theo}

\begin{theo} \label{T8.6} 
If $\lambda = \frac{1-n}{2}$ and $\mu = \frac{1+n}{2}$ for some $n \in \nit_0$, then for
each $u$, one has a short exact sequence
$$0 \rightarrow H^u(s\ell_2, \D_{\lambda\mu}^{n-1}(\rit)) \stackrel{i_\sharp}{\rightarrow}
H^u(s\ell_2, \D_{\lambda\mu}^{n}(\rit)) \stackrel{\sigma_\sharp}{\rightarrow} H^u(s\ell_2,
C_\infty(\rit,S_n^n \rit)) \rightarrow 0$$
where the kernel (resp. the quotient) is isomorphic to $\rit$ for $u=1,2$ (resp. $0,1$)
and is vanishing otherwise. Moreover, the inclusion of $\D_{\lambda\mu}^{n}(\rit)$ into
$\D_{\lambda\mu}(\rit)$ induces an isomorphism in cohomology.
\end{theo}

\begin{rema} \label{R8.7}{\rm
1) One sees in particular that $H^3(s\ell_2,\D_{\lambda\mu}(\rit))=0$ for all $\lambda,\mu
\in \rit$.

2) For $(\lambda,\mu) = (\frac{1-n}{2},\frac{1+n}{2})$, $n \in \nit_0$,
$H^1(s\ell_2,\D_{\lambda\mu}(\rit))$ is spanned by
$$\tau_n' \left(f \frac{d}{dt}\right) : g \mod{dt}^\lambda \rightarrow \frac{df}{dt}\;
\frac{d^ng}{dt^n}\;\mod{dt}^\mu$$
and
$$\gamma_n' \left(f \frac{d}{dt}\right) : g \mod{dt}^\lambda \rightarrow \frac{d^2f}{dt^2}\;
\frac{d^{n-1}g}{dt^{n-1}}\;\mod{dt}^\mu.$$
This follows from Theorem \ref{T8.6}. The $\gamma_n'$'s are exactly the cocycles used in
\cite{3} in the study of the $s\ell_2$-structure of $\D_{\lambda\mu}(\rit)$.

3) If $p=q+n$ and $\delta = \frac{p+q+n}{2}$, through the isomorphism between
$\D(\S_\delta^p(\rit),\S_\delta^q(\rit))$ and $\D_{\frac{1-n}{2},\frac{1+n}{2}}(\rit)$,
$\tau_n$ and $\gamma_n$ become respectively $\tau_n'$ and $\gamma_n'$.}
\end{rema}

\subsection*{Acknowledgements}

\indent\indent We would like to thank very much M. De Wilde, P. Mathonet and V. Ovsienko 
for helpful discussions.

\begin{tabular}{l}
P. LECOMTE \\
Institut de MathŽmatique \\
UniversitŽ de Lige \\
Grande Traverse 12, B37 \\
B-4000 Lige (Belgium) \\ \\
e-mail : plecomte@ulg.ac.be
\end{tabular}

\end{document}